\documentclass[12pt]{article} 
\usepackage{amsmath,amsthm, amssymb} 
\usepackage{amssymb,latexsym}

\newtheorem{theorem}{Theorem}

\newtheorem{lemma}{Lemma}

\begin{document}

\baselineskip=17pt 

\title{\bf The binary Goldbach problem \\ with arithmetic weights attached \\ to one of the variables}

\author{\bf D. I. Tolev \footnote{Supported by Sofia University Grant 221}}

\date{}
\maketitle

\section{Introduction and statement of the results.}

Suppose that $N$ is a sufficiently large integer and denote
\[
J(n) = \sum_{p_1 + p_2 = n} \log p_1 \log p_2 .
\]
(From this place the letter $p$, with or without subscripts, is reserved for primes.)
It is expected that if $n$ is a large even integer then $J(n) \sim c_0 \lambda (n) n $, where
\begin{equation} \label{+++1}
\lambda(k) = \prod_{\substack{p \mid k \\ p > 2}} \frac{p-1}{p-2} , \qquad
 c_0  =   2 \prod_{p>2} \left( 1 - \frac{1}{(p-1)^2} \right) .
\end{equation}
This conjecture has not been proved so far, but using the Hardy--Littlewood circle method 
and Vinogradov's method for estimating exponential sums over primes
(see, for example,  Vaughan~\cite{Vaug}, Ch.~2), one can find that
\begin{equation} \label{+++2}
\sum_{\substack{n \le N \\ 2 \mid n}} \left| J(n) - c_0 \lambda (n) n \right|
\ll N^2 \mathcal L^{-A} ,
\end{equation}
where $A>0$ is an arbitrarily large constant and $\mathcal L = \log N$.

Let $r(k)$ be the number of solutions of the equation $x_1^2 + x_2^2 = k$ in integers $x_1, x_2$.
One of the classical problems in prime number theory is the Hardy--Littlewood problem
concerning the representation of large integers as a sum of two squares and a prime.
It was solved by Linnik (see \cite{Linnik}) and related problems have been studied
by Linnik, Hooley and other mathematicians.
For more information we refer the reader to Hooley's book \cite{Hooley},~Ch.5.
In particular, one can show that
\begin{equation} \label{+++3}
   \sum_{p \le N} r(p - 1) = 
   \pi N \mathcal L^{-1}
    \prod_{p > 2} \left( 1 + \frac{\chi(p)}{p(p-1)} \right)
    +    O \left( N \mathcal L^{-1-\theta_0} \left( \log \mathcal L \right)^5 \right) ,
\end{equation}
where $ \chi(k)$ is the non-principal character modulo 4 and
\begin{equation} \label{+++4}
\theta_0 = \frac{1}{2} - \frac{1}{4} e \log 2 = 0.0029\dots .
\end{equation}

Let $\tau (k)$ be the number of positive divisors of $k$. 
Linnik~\cite{Linnik} 
(see also Halberstam and Richert~\cite{H-R},~Ch.~3.5.)
solved the Titchmarsh divisor problem and proved
that
\begin{equation} \label{+++5}
   \sum_{p \le N} \tau(p - 1)  = c_0 N  + 
     O \left( N \mathcal L^{-1} \log \mathcal L \right) , \qquad
     c_0 = \prod_p \left( 1 + \frac{1}{p (p-1)}\right) .
\end{equation}
We note that sharper versions of \eqref{+++3} and \eqref{+++5} are known at present
(see Bredihin~\cite{Bredihin}, 
Bombieri, Friedlander and Iwaniec~\cite{BoFrIw}) and Fouvry~\cite{Fouvry}.

In this paper we state two theorems which are, in some sense, combinations of
\eqref{+++2}, \eqref{+++3} and respectively \eqref{+++2}, \eqref{+++5}.
Denote 
\begin{equation} \label{+++5.5}
\mathcal R(n) = \sum_{p_1 + p_2 = n} r(p_1 - 1) \log p_1 \log p_2 .
\end{equation}
After certain formal calculations one may conjecture that for any 
sufficiently large even $n$ the quantity
$\mathcal R(n)$ is asymptotically equal to 
\begin{equation} \label{+++6}
  \mathcal M_{\mathcal R} (n) 
   = \pi c_0 n 
    \prod_{p \mid n-1} \left( 1 - \frac{\chi(p)}{p} \right)
  \prod_{\substack{p \mid n \\ p > 2 }} \left( 1 + \frac{p + \chi(p) }{p(p-2)} \right)
  \prod_{p \nmid n(n-1)} \left( 1 + \frac{2 \chi(p)}{p(p-2)} \right) .
\end{equation}
Our first result is the following:
\begin{theorem} \label{T1}
Suppose that $\theta_0$ is the constant defined by \eqref{+++4}. Then
we have
\begin{equation} \label{+++7}
\sum_{\substack{ n \le N \\ 2 \mid n } } 
  \left|
    \mathcal R(n) -  \mathcal M_{\mathcal R} (n)  
  \right|
   \ll N^2 \mathcal L^{- \theta_0}  \left( \log \mathcal L \right)^6 .
\end{equation}
\end{theorem}
It is clear that $  n \left( \log \log (10 n) \right)^{-2} \ll \mathcal M_{\mathcal R} (n) \ll n \left( \log \log (10 n) \right)^2$. Also, from \eqref{+++7} it follows that for any positive constant $\theta <  \theta_0$ 
the number of even $n \le N$ for which 
$\left|  \mathcal R(n) -  \mathcal M_{\mathcal R} (n)  \right| > N \mathcal L^{-\theta} $
is $ O \left(N \mathcal L^{- (\theta_0 - \theta) }  (\log \mathcal L )^6 \right)$. 
So, in other words, $\mathcal R(n)$ is close to 
$ \mathcal M_{\mathcal R} (n) $ for almost all even $n$.

Theorem~\ref{T1} is related to a recent result of K. Matom\"aki~\cite{Mato}. It is shown in \cite{Mato} that 
the number of integers $n \le N$ satisfying $n \equiv 0 \; \text{or} \; 4 \pmod{6} $ and
that cannot be represented as a sum of two primes, one of which of the form
$k^2 + l^2 +1$, is $O \left( N \mathcal L^{-A} \right)$, where $A$ is an arbitrarily large constant. 
So Matom\"aki's estimate for the cardinality of this exceptional set is stronger then ours, but her method
does not provide so sharp information about the number of such representations.

Our second result is concerning the quantity
\[
\mathcal T(n)
 = \sum_{p_1 + p_2 = n} \tau(p_1 - 1) \log p_1 \log p_2 .
\]
Again, after certain formal calculations, one may conclude that
$\mathcal T(n)$ should be asymptotically equal to
\[
\mathcal M_{\mathcal T}(n) = c_0 n \log n 
    \prod_{p \mid n-1} \left( 1 - \frac{1}{p} \right)
  \prod_{\substack{p \mid n \\ p > 2 }} \left( 1 + \frac{p + 1 }{p(p-2)} \right)
  \prod_{p \nmid n(n-1)} \left( 1 + \frac{2}{p(p-2)} \right) .
\]
We can establish:
\begin{theorem} \label{T2}
The following estimate holds
\[
 \sum_{\substack{ n \le N \\ 2 \mid n } } 
  \left|
    \mathcal T(n) - \mathcal M_{\mathcal T} (n)  
  \right|
   \ll N^2 \left( \log \mathcal L \right)^3 .
\]
\end{theorem}
We note that 
$ n \, \log n \, \left( \log \log (10 n) \right)^{-2} \ll \mathcal M_{\mathcal T} (n) \ll n \, \log n \,
 \left( \log \log (10 n) \right)^2 $, 
so the quantity $\mathcal T(n)$ is close to $ \mathcal M_{\mathcal T} (n) $
for almost all even $n$.

We prove only Theorem~\ref{T1}. The proof of Theorem~\ref{T2} is similar and simpler.

\section{Some lemmas.}

Suppose that $n \le N$ and let $k$ and $l$ be integers with
$(k,l)=1$ (as usual, $(k,l)$ stands for the greatest common factor
of $k$ and $l$). 
Let  $\mathcal I $ be the set of all subintervals of the interval $[1, N]$ and let $I \in \mathcal I$.
We denote
\begin{align} 
J_{k, l}(n; I) 
  & = 
  \sum_{\substack{p_1 + p_2 = n \\ p_1 \equiv l \pmod{k} \\ p_1 \in I}} 
\log p_1 \log p_2 , 
\qquad
J_{k, l}(n) = J_{k, l} (n; [1, N]) ;
  \label{+++8} \\
\mathfrak S_{k, l}(n)
  & = 
 \begin{cases}
   c_0 
   \lambda(nk)
       \qquad
      &    \text{if} \;\; (k, n-l) = 1 \; \; \text{and} \;\;2 \mid n , \\
       0 
       &  \text{otherwise} ;
      \end{cases} 
          \label{+++9} \\
  \Phi(n; I) 
   & = 
   \sum_{\substack{m_1 + m_2 = n \\ m_1 \in I}} 1 .
   \label{+++9.5}
\end{align}

Our first lemma states that the expected formula for $J_{k, l}(n; I) $ is true on average with respect to
$k \le \sqrt{N} \mathcal L^{-B}$ and $n \le N$ and uniformly for $l$ and $I$. More precisely, we have
\begin{lemma} \label{L1}
For any constant $A>0$ there exist $B=B(A) > 0$ such that
\[
\sum_{k \le \sqrt{N} \mathcal L^{-B} } \max_{(l, k)=1} \max_{I \in \mathcal I}
  \sum_{ n \le N } \left|
   J_{k, l}(n; I) - \frac{\mathfrak S_{k, l}(n)}{\varphi(k)} \Phi(n; I)
     \right|
     \ll N^2 \mathcal L^{-A} .
\]
\end{lemma}

This lemma is very similar to results of Mikawa~\cite{Mikawa} and Laporta~\cite{Laporta}. 
These authors study the equation $p_1 - p_2 = n$ and without the condition $p_1 \in I$. 
However inspecting the arguments presented in  \cite{Laporta}, 
the reader will readily see that the proof of Lemma~\ref{L1} can be obtained is the same manner.

The next lemma is an immediate consequence from a classical sieve theory result
(see \cite{H-R},~Ch.~2, Th.~2.4).

\begin{lemma} \label{L2}
Suppose that $h$ is an integer such that $1 \le |h| \le N$. Then the number of 
solutions of the equation $p_1 - p_2 = h$ in primes $p_1, p_2 \le N$ is
$O \left( N \mathcal L^{-2} \log \mathcal L \right)$, where the constant in the Landau symbol is absolute.
\end{lemma}

The next two lemmas are due to C.Hooley and play an essential role in the proof of \eqref{+++3}, 
as well as in the solutions of other related problems. 
\begin{lemma} \label{L3}
Suppose that $\omega > 0 $ is a constant and let $F_{\omega}(N)$ be the number of primes $p \le N$ such that
$p-1$ has a divisor lying between $ \sqrt{N} \mathcal L^{-\omega} $ and $ \sqrt{N} \mathcal L^{\omega} $.
Then we have
\[
F_{\omega}(N) \ll N \mathcal L^{-1 - 2 \theta_0} \left( \log \mathcal L \right)^3 ,
\]
where $\theta_0$ is defined by \eqref{+++4} and where
the constant in the Vinogradov symbol depends only on $\omega$.
\end{lemma}

\begin{lemma} \label{L4}
Suppose that $\omega > 0 $ is a constant. Then we have
\[
\sum_{p \le N} \left| 
   \sum_{\substack{ d \mid p - 1 \\\sqrt{N} \mathcal L^{-\omega} < d <  \sqrt{N} \mathcal L^{\omega} } }
   \chi(d)
\right|^2
   \ll N \mathcal L^{-1} \left( \log \mathcal L \right)^7 ,
\]
where the constant in the Vinogradov symbol depends only on $\omega$.
\end{lemma}

The proofs of very similar results (with $\omega = 48$ and with the condition $d \mid N - p$ 
rather than $d \mid p-1$) are available in \cite{Hooley},~Ch.5 and the reader will easily see that 
the method used there yields also the validity of Lemmas~\ref{L3} and~\ref{L4}.

\section{Proof of Theorem~\ref{T1}.}

\subsection{Beginning.}

Denote by $\mathcal E$ the sum on the left-hand side of \eqref{+++7} and put
\begin{equation} \label{+++10}
  D = \sqrt{N} \mathcal L^{-1 - B(1)} ,
\end{equation}
where $B(A)$ is specified in Lemma~\ref{L1}.
Using \eqref{+++5.5} and the well-known identity
$
r(m) = 4 \sum_{d \mid m} \chi(d)
$
we find
\begin{equation} \label{+++11}
\mathcal R(n) = 4 
  \sum_{p_1 + p_2 = n} \left( 
   \sum_{d \mid p_1 -1} \chi(d)   \right)
  \log p_1 \log p_2
   = 4 \left( S_1(n) + S_2(n) + S_3(n) \right) ,
\end{equation}
where
\begin{align}
 S_1(n)
  &= 
   \sum_{p_1 + p_2 = n} \left( 
   \sum_{\substack{d \mid p_1 -1 \\ d \le D}} \chi(d)   \right)
  \log p_1 \log p_2
  \label{+++12} \\
 S_2(n)
  &= 
   \sum_{p_1 + p_2 = n} \left( 
   \sum_{\substack{d \mid p_1 -1 \\ D < d < N/D}} \chi(d)   \right)
  \log p_1 \log p_2
  \label{+++13} \\
  S_3(n)
  &= 
   \sum_{p_1 + p_2 = n} \left( 
   \sum_{\substack{d \mid p_1 -1 \\ d \ge N/D}} \chi(d)   \right)
  \log p_1 \log p_2
  \label{+++14}
\end{align}

Therefore from \eqref{+++7} and \eqref{+++11} it follows
\begin{equation} \label{+++15}
\mathcal E \ll 
\mathcal E_1 + \mathcal E_2 + \mathcal E_3 ,
\end{equation}
where 
\begin{equation} \label{+++16}
\mathcal E_1 
   = 
  \sum_{\substack{n \le N \\ 2 \mid n}} 
  \left|
    4 S_1(n)   - \mathcal M_{\mathcal R}(n) 
  \right| ;
     \qquad
   \mathcal E_j
   =
   \sum_{\substack{n \le N \\ 2 \mid n }} 
  \left|   S_j(n)   \right|  , \qquad j = 2,3.
\end{equation}

\subsection{The estimation of $\mathcal E_1$.}

Using \eqref{+++8}, \eqref{+++9.5}, \eqref{+++12} and bearing in mind Lemma~\ref{L1} we find
\[
S_1(n) = \sum_{d \le D} \chi(d) J_{d, 1}(n) = (n-1) S_1^{\prime}(n) + S_1^*(n) ,
\]  
where
\begin{align} 
 S_1^{\prime}(n) 
  & 
  = \sum_{d \le D} \chi(d) 
 \frac{\mathfrak S_{d, 1}(n)}{\varphi(d)} , 
  \label{+++17} \\
  S_1^*(n) 
   & = 
   \sum_{d \le D} \chi(d) 
 \left( J_{d, 1}(n) - (n-1) \frac{\mathfrak S_{d, 1}(n)}{\varphi(d)}
  \right) .
  \label{+++17.5}
\end{align}
Hence
\begin{equation} \label{+++18}
\mathcal E_1 \ll \mathcal E_1^{\prime} +  \mathcal E_1^* ,
\end{equation}
where
\begin{equation} \label{+++19}
\mathcal E_1^{\prime} = \sum_{\substack{n \le N \\ 2 \mid n}} 
     \left| 
   4 (n-1) S_1^{\prime}(n) -  \mathcal M_{\mathcal R}(n)   \right| , \qquad
   \mathcal E_1^* =  \sum_{\substack{n \le N \\ 2 \mid n}} \left| 
   S_1^*(n)   \right| .
\end{equation}

By  \eqref{+++10}, \eqref{+++17.5}, \eqref{+++19} and Lemma~\ref{L1} it follows that
\begin{equation} \label{+++20}
 \mathcal E_1^* \ll N^2 \mathcal L^{-1} .
\end{equation}

Consider $\mathcal E_1^{\prime}$.
From \eqref{+++1}, \eqref{+++9} and \eqref{+++17} we find
\begin{equation} \label{+++21}
S_1^{\prime}(n) = c_0
   \sum_{\substack{d \le D \\ (d, n-1)=1}} \frac{\chi(d)}{\varphi(d)} 
   \lambda (nd)
       = c_0 \lambda(n)
   \sum_{\substack{d \le D \\ (d, n-1)=1}} f_n(d) ,
\end{equation}
where
\begin{equation} \label{+++22}
 f_n(d) = \frac{\chi(d)}{\varphi(d)} 
   \frac{\lambda (d)}{\lambda((n, d))} .
\end{equation}
Obviously the function $f_n(d)$ is multiplicative with respect to $d$ and 
\begin{equation} \label{+++23}
f_n(d) \ll d^{-1} \left( \log \log (10 d) \right)^2 
\end{equation}
uniformly with respect to $n$. 
To evaluate the sum in right-hand side of \eqref{+++21} we consider the function
\[
F_n(s) = \sum_{\substack{d=1 \\ (d, n-1) = 1 }}^{\infty} f_n(d) d^{-s} .
\]
It is analytic in the half-plane $ Re \, (s) > 0 $ and we may represent it as an Euler product:
\[
F_n(s) = \prod_{p \nmid n-1} T_n(p, s) , \qquad T_n(p, s) = 1 + \sum_{l=1}^{\infty}  f_n(p^l) p^{-ls} .
\]

From \eqref{+++1} and \eqref{+++22} we easily find
\[
f_n(p^l) =
 \begin{cases}
   \chi(p)^l \; p^{1-l} \; (p-1)^{-1}
    & \text{if} \; \; p \mid n , \\
   \chi(p)^l \; p^{1-l} \; (p-2)^{-1}
    & \text{if} \; \; p \nmid n  ;
 \end{cases}
\]
and respectively
\[
T_n(p, s) =
  \left( 1 - \frac{\chi(p)}{ p^{s+1}} \right)^{-1} 
  T_n^*(p, s) ,
\]
where
\[
T_n^*(p, s) = 
  \begin{cases}
       1 +  \chi(p)  p^{-s-1} (p-1)^{-1} 
     & \quad \text{if} \quad  p \mid n , \\
       1 + 2 \chi(p)  p^{-s-1} (p-2)^{-1} 
     & \quad \text{if} \quad  p \nmid n .
     \end{cases}
\]

Therefore
\begin{equation} \label{+++24}
F_n(s) = L(s+1, \chi) H_n(s)
\end{equation}
where $L(s, \chi) $ is the Dirichlet $L$-function corresponding to the character $\chi$ and
\begin{equation} \label{+++25}
 H_n(s) =
  \prod_{p \mid n-1}
  \left( 1 - \frac{\chi(p)}{ p^{s+1}} \right)
  \;
  \prod_{p \mid n}
  \left( 1 + \frac{\chi(p)}{ p^{s+1} (p-1)} \right)
  \;
  \prod_{p \nmid n(n-1)}
  \left( 1 + \frac{2 \chi(p)}{ p^{s+1} (p-2)} \right) .
\end{equation}

From \eqref{+++24}, \eqref{+++25} we see that $F_n(s)$ has an analytic continuation to the half-plane
$Re \, (s) > -1$. 
It is clear that $H_n(s) \ll n^{\varepsilon}$ for $|Re \, (s)| \ge -1/2$
(here and later $\varepsilon$ is an arbitrarily small positive number).
Also, it is well-known that in the same region 
we have $L(s+1, \chi) \ll 1 + |Im \, (s)|^{1/6}$. Hence
\begin{equation} \label{+++26}
F_n(s) \ll N^{\varepsilon} \, T^{1/6} \qquad \text{if} \qquad Re \, (s) \ge -1/2 , \quad | Im \, (s) | \le T 
\end{equation}
for any $T > 1 $.
We apply Perron's formula (see, for example \cite{Tenen},~Ch.~II.2) to find
\begin{equation} \label{+++27}
\sum_{\substack{d \le D \\ (d, n-1)=1}} f_n(d) = \frac{1}{2 \pi i} 
\int_{\varkappa - i T}^{\varkappa + i T}
  F_n(s) \frac{D^s}{s} d s + 
  O \left( \sum_{d=1}^{\infty} 
\frac{D^{\varkappa} \, |f_n(d)|}{  d^{\varkappa} \left(1 + T \left| \log \frac{D}{d} \right| \right)}  \right) 
\end{equation}
with $\varkappa = 1/10$ and $T = N^{3/4}$.
Using \eqref{+++10} and \eqref{+++23} one can easily verify that
the remainder term in \eqref{+++27} is $O \left( N^{-1/20}\right) $.
To evaluate the integral in \eqref{+++27}
we apply Cauchy's theorem. The residue of the integrand at $ s=0 $ equals
\begin{equation} \label{+++28}
F_n(0) =
\frac{\pi}{4} 
  \prod_{p \mid n-1} \left( 1 - \frac{\chi(p)}{p} \right) 
  \prod_{p \mid n} \left( 1 + \frac{\chi(p)}{p(p-1)} \right) 
  \prod_{p \nmid n(n-1) } \left( 1 + \frac{2 \chi(p)}{p(p-2)} \right) .
\end{equation}
Hence the main term in the right-hand side of \eqref{+++27} is equal to
\begin{equation} \label{+++29}
F_n(0)
 + \frac{1}{2 \pi i}
   \left(
    \int_{\varkappa - i T}^{-1/2 - i T} + \int_{-1/2 - i T}^{-1/2 + i T} + \int_{-1/2 + i T}^{\varkappa + i T}     
   \right)
     F_n(s) \frac{D^s}{s} d s .
\end{equation}
Using \eqref{+++26} one can easily find that the contribution of the integrals in \eqref{+++29}
is $O \left( N^{-1/20}\right)$. 
Therefore 
\begin{equation} \label{+++30}
\sum_{\substack{d \le D \\ (d, n-1)=1}} f_n(d) = F_n(0) + O \left( N^{-1/20}\right) . 
\end{equation}
From \eqref{+++1}, \eqref{+++6}, \eqref{+++19}, \eqref{+++21}, \eqref{+++28} and \eqref{+++30} it follows that
\[
\mathcal E_1^{\prime} \ll N^2 \mathcal L^{-1} .
\]
Hence, using \eqref{+++18} and \eqref{+++20} we get
\begin{equation} \label{+++31}
  \mathcal E_1 \ll N^2 \mathcal L^{-1} .
\end{equation}

\subsection{The estimation of $\mathcal E_2$.}

Clearly,  from \eqref{+++16} and Cauchy's inequality it follows that
\begin{equation} \label{+++32}
\mathcal E_2 \ll N^{1/2} \left( \sum_{n \le N} \left| S_2(n) \right|^2 \right)^{1/2} 
 = N^{1/2} \left( \mathcal E_2^{\prime} \right)^{1/2} ,
\end{equation}
say. 
Using \eqref{+++13} we find
\begin{align}
\mathcal E_2^{\prime} 
  & = 
\sum_{n \le N} \sum_{D < d, t < N/D}
 \chi(d) \chi(t) 
 \sum_{\substack{p_1 + p_2 = n \\ p_1 \equiv 1  \pmod{d} }} \log p_1 \log p_2
 \sum_{\substack{p_3 + p_4 = n \\ p_3 \equiv 1  \pmod{t} }} \log p_3 \log p_4
  \notag \\
  & \notag \\
  & =
  \sum_{p_1 + p_2 = p_3 + p_4 \le N} \log p_1 \log p_2 \log p_3 \log p_4
 \sum_{\substack{D < d, t < N/D  \\ d \mid p_1 - 1 , \;\; t \mid p_3 - 1}}
 \chi(d) \chi(t) 
 \notag \\
 & \notag \\
 & \ll
 \mathcal L^4 \,
    \mathcal E_2^{\prime\prime} + N^{2 + \varepsilon} ,
 \label{+++33}
\end{align} 
where
\[
  \mathcal E_2^{\prime\prime}
   =
 \sum_{\substack{p_1 + p_2 = p_3 + p_4 \\ p_1, p_2, p_3, p_4 \le N \\ p_1 \not= p_3}}
  \left|
  \sum_{\substack{D < d < N/D  \\ d \mid p_1 - 1 }}
 \chi(d) 
  \right|
  \; 
   \left|
  \sum_{\substack{D < t < N/D  \\ t \mid p_3 - 1 }}
 \chi(t) 
  \right| .
\]

Denote by $\mathcal F$ the set of primes $p \le N$ such that $p-1$ has a divisor
lying between $D$ and $N/D$.
Using the inequality $uv \le u^2 + v^2$
and taking into account the symmetry with respect to $d$ and $t$
we get
\begin{align}
\mathcal E_2^{\prime\prime}
   & \ll 
  \sum_{\substack{p_1 + p_2 = p_3 + p_4 \\ p_1, p_2, p_4 \le N \\ p_1 \not= p_3 , \; p_3 \in \mathcal F }}
\left|
  \sum_{\substack{D < d < N/D  \\ d \mid p_1 - 1 }}
 \chi(d) 
  \right|^2
  \notag \\
  & =
    \sum_{p_1 \le N}
  \left|
     \sum_{\substack{D < d < N/D  \\ d \mid p_1 - 1 }}
 \chi(d) 
  \right|^2
\sum_{\substack{ p_3 \in \mathcal F \\ p_3 \not= p_1 }}  \;
  \sum_{\substack{p_2, p_4 \le N \\ p_4 - p_2 = p_1 - p_3}} 1 .
   \label{+++33.5}
\end{align}  

Applying Lemmas~\ref{L2} and \ref{L3} we find
\begin{equation} \label{+++33.6}
\sum_{\substack{p_3 \in \mathcal F \\ p_3 \not= p_1 }}
  \sum_{\substack{p_2, p_4 \le N \\ p_4 - p_2 = p_1 - p_3}} 1 
    \ll
    N \mathcal L^{-2} ( \log \mathcal L ) \sum_{p \in \mathcal F} 1 \ll
      N^2 \mathcal L^{-3- 2 \theta_0} ( \log \mathcal L )^4 
\end{equation}
and then using \eqref{+++33.5}, \eqref{+++33.6} and 
Lemma~\ref{L4} we get
\begin{equation} \label{+++34}
\mathcal E_2^{\prime\prime}
 \ll
   N^2
   \mathcal L^{-3 - 2 \theta_0} (\log \mathcal L)^4
   \sum_{p \le N}
\left|
  \sum_{\substack{D < d < N/D  \\ d \mid p - 1 }}
 \chi(d) 
  \right|^2 
  \ll N^3
   \mathcal L^{ -4 - 2 \theta_0} (\log \mathcal L )^{11} .
\end{equation}

From \eqref{+++32}, \eqref{+++33} and \eqref{+++34} we conclude that 
\begin{equation} \label{+++35}
\mathcal E_2 \ll N^2 \mathcal L^{-\theta_0 }  \left( \log \mathcal L \right)^6 .
\end{equation}

\subsection{The estimation of $\mathcal E_3$}

From \eqref{+++14} it follows that
\begin{align}
S_3(n) 
  & = 
  \sum_{p_1 + p_2 = n} \log p_1 \log p_2 
     \sum_{\substack{ m \mid p_1 - 1 \\ \frac{p_1-1}{m} \ge \frac{N}{D} }} \chi \left( \frac{p_1 - 1 }{m} \right)
     \notag \\
  & = 
  \sum_{p_1 + p_2 = n} \log p_1 \log p_2 
    \sum_{j = \pm 1} \chi(j)
     \sum_{\substack{m \le \frac{(p_1-1) D}{N} , \; \; 2 \mid m \\ p_1 \equiv 1 + jm  \pmod{4m} }} 1 .
     \notag 
\end{align}         
We change the order of summation and use \eqref{+++8} to find
\[
   S_3(n) = 
    \sum_{\substack{m \le D \\ 2 \mid m }} \sum_{j = \pm 1} \chi(j)
    J_{4m, 1+jm}(n, I_m) ,
\]
where $I_m$ denotes the interval $\left[1 + m N/D , N \right]$. Having in mind Lemma~\ref{L1}
we write
\begin{equation} \label{+++35.5}
S_3(n) = S_3^{\prime}(n) + S_3^*(n) ,
\end{equation}
where
\begin{align}
S_3^{\prime}(n) 
  & =   
  \sum_{\substack{ m \le D \\ 2 \mid m }} \sum_{j = \pm 1} \chi(j) \frac{\mathfrak S_{4m, 1 + jm}(n)}{\varphi (4m)}
  \Phi(n, I_m) ,
  \notag \\
S_3^*(n) 
  & =
  \sum_{\substack{ m \le D \\ 2 \mid m }} \sum_{j = \pm 1} \chi(j)
    \left(
    J_{4m, 1+jm}(n, I_m) - \frac{\mathfrak S_{4m, 1 + jm}(n)}{\varphi (4m)}   \Phi(n, I_m) 
        \right) .
  \label{+++35.6}
\end{align}

Since $2 \mid n$ it follows from \eqref{+++9} that 
\[
  \mathfrak S_{4m, 1 + jm}(n)  = 
   \begin{cases}  
      c_0 \lambda(4mn) & \text{if} \;\; (4m, n-1-jm)=1 , \\ 
      0                & \text{otherwise} . 
      \end{cases}
\]
However the condition
$(4m, n-1-jm) = 1$ is independent of $j$ 
(from the set $\{  1, -1 \} $) and therefore
$\mathfrak S_{4m, 1 + jm}(n) $ is independent of $j$ too. This means that
\[
S_3^{\prime}(n) =  0 .
\]
Hence, using \eqref{+++10}, \eqref{+++16}, \eqref{+++35.5}, \eqref{+++35.6} and Lemma~\ref{L1} we find
\begin{align}
\mathcal E_3 
  & \ll 
  \sum_{n \le N} \left| S_3^*(n) \right|
    \notag \\
    & \notag \\
  &  \ll
 \sum_{\substack{m \le D \\ 2 \mid m }} \sum_{j = \pm 1} \sum_{n \le N}
    \left|
    J_{4m, 1+jm}(n, I_m) - \frac{\mathfrak S_{4m, 1 + jm}(n)}{\varphi (4m)}   \Phi(n, I_m) 
        \right|
     \notag \\
     & \notag \\
  & \ll
  \sum_{k \le 4D}  \max_{(l, k)=1} \max_{I \in \mathcal I} \sum_{n \le N}
   \left|
    J_{k, l}(n, I) - \frac{\mathfrak S_{k, l}(n)}{\varphi (k)}   \Phi(n, I) 
        \right|
   \notag \\
   & \notag \\
   & \ll
   N^2 \mathcal L^{-1} .
        \label{+++36}
\end{align}

The estimate \eqref{+++7} follows from \eqref{+++15},
\eqref{+++31}, \eqref{+++35} and \eqref{+++36}, so the theorem is proved.

\bigskip
\bigskip

\vbox{
\hbox{Faculty of Mathematics and Informatics}
\hbox{Sofia University ``St. Kl. Ohridsky''}
\hbox{5 J.Bourchier, 1164 Sofia, Bulgaria}
\hbox{ }
\hbox{Email: dtolev@fmi.uni-sofia.bg}}

\end{document}